\documentclass[a4paper,12pt]{article}

\usepackage[utf8]{inputenc}
\usepackage{amsmath, amsfonts, amssymb, amsthm}
\usepackage{mathrsfs}
\usepackage{fullpage}

\usepackage{indentfirst}

\usepackage{bbm}

\numberwithin{equation}{section}

\newtheorem{theorem}{Theorem}
\newtheorem{fact}{Fact}
\newtheorem*{lemma*}{Lemma}
\newtheorem*{defin}{Definition}

\newcommand{\Z}{\mathbb{Z}}
\newcommand{\sign}{\mathrm{sign}}
\newcommand{\wh}{\widehat}
\newcommand{\F}{\mathcal{F}}
\newcommand{\E}{\mathbb{E}}
\newcommand{\eps}{\varepsilon}
\newcommand{\LPR}{\mathrm{LPR}}
\newcommand{\Rad}{\mathrm{Rad}}
\newcommand{\wt}{\widetilde}
\newcommand{\kap}{\varkappa}
\newcommand{\Sq}{\mathcal{S}}

\title{Some notes on the vector-valued extension of Littlewood--Paley--Rubio de Francia inequality for Walsh functions}

\author{Anton Tselishchev\thanks{This research was supported by the Russian Science Foundation (grant No.~18-11-00053).}}

\date{}

\begin{document}
	\maketitle
	
	\begin{abstract}
			J. L. Rubio de Francia proved the one-sided Littlewood--Paley inequality for arbitrary intervals in $L^p$, $2\le p<\infty$ and later N. N. Osipov proved the similar inequality for Walsh functions. In this paper we investigate some properties of Banach spaces $X$ such that the latter inequality holds for $X$-valued functions.
	\end{abstract}
	
	\section{Introduction}
	
	Let $\{I_m\}$ be a sequence of disjoint intervals in $\Z$. In the paper \cite{RdF} Rubio de Francia showed that the following inequality holds for any function $f\in L^p(\mathbb{T})$, $2\le p<\infty$:
	\begin{equation}
	\Big\| \Big( \sum_{m} |(\widehat{f}\mathbbm{1}_{I_m})^{\vee}|^2 \Big)^{1/2} \Big\|_{L^p}\lesssim \|f\|_{L^p}.
	\label{RdF_exp}
	\end{equation}
	We use the notation ``$\lesssim$" to indicate that the left hand side does not
	exceed some positive constant times the right hand side. Here this constant does not depend on $f$ or the intervals $I_m$. It is worth noting that Rubio de Francia worked with functions defined on $\mathbb{R}$ but it is easy to transfer the results of the paper \cite{RdF} to the circle $\mathbb{T}$.
	
	Later, in the paper \cite{Osip} N. Osipov proved the same inequality (in the ``dual" form) in the context of Walsh functions. Since we will use Walsh functions in what follows, we remind the reader the basic notions concerning them.
	
	The Rademacher functions are defined as $r_k(x)=\sign\sin 2^k\pi x$. For $n=2^{k_1}+2^{k_2}+\ldots+2^{k_m}$ where $k_1>k_2>\ldots>k_m\ge 0$ we define the $n$th Walsh function as $w_n=r_{k_1+1}r_{k_2+1}\ldots r_{k_m+1}$. It is well-known that Walsh functions form an orthonormal basis in $L^2[0,1]$. For any function $f$ defined on $[0,1]$ we denote by $\wh{f}$ the sequence of its Walsh coefficients, that is, $\wh{f}(n)=(f,w_n)=\int fw_n$.
	
	Suppose that 
	$$
	n_1=\sum_{k=0}^m \alpha_k 2^k \quad \text{and}\quad n_2=\sum_{j=0}^m \beta_j 2^j
	$$
	are two integers and $\alpha_k$ and $\beta_j$ are their binary digits. We are going to use the following notation:
	$$
	n_1\dotplus n_2 = \sum_{k=0}^m ((\alpha_k+\beta_k)\mathrm{mod}\, 2) 2^k.
	$$
	Then it is easy to see that the following formula holds true:
	$$
	w_{n_1}w_{n_2}=w_{n_1\dotplus n_2}.
	$$

	The Walsh functions are closely related to dyadic martingales. We denote by $\F_k$ the $\sigma$-algebra generated by dyadic intervals of length $2^{-k}$. It is easy to see that the conditional expectation of a function $f$ with respect to $\F_k$ may be written as
	$$
	\E(f| \F_k)=\sum_{n=0}^{2^k-1}(f,w_n)w_n.
	$$
	We will use the notation $\E_k f$ instead of $\E(f|\F_k)$ for simplicity. If we denote by $\delta_k$ the interval $[2^{k-1}, 2^k-1]$ in $\Z_+$, then the martingale differences may be written in the following way:
	$$
	\Delta_k f=\E_k f-\E_{k-1}f=\sum_{n\in\delta_k} (f,w_n)w_n, \qquad k\ge 1.
	$$
	We also put $\delta_0=\{0\}$ and $\Delta_0 f=(f,w_0)w_0$.
	
	For any set $A\subset\Z_+$ we define the following orthogonal projection on $L^2[0,1]$:
	$$
	P_A f=\sum_{n\in A}(f,w_n)w_n=(\chi_A\wh{f})^{\vee}.
	$$
	
	If $\{I_s\}_s$ is a collection of pairwise disjoint intervals in $\Z_+$, then the following analogue of Rubio de Francia's inequality \eqref{RdF_exp} holds for any function $f\in L^p[0,1]$, $2\le p<\infty$:
	\begin{equation}
	\Big\|\Big(\sum_s |P_{I_s}f|^2\Big)^{1/2}\Big\|_{L^p}\lesssim \|f\|_{L^p}.
	\label{RdF_walsh}
	\end{equation}
As we mentioned before, the proof of this inequality (in the dual form) may be found in the paper \cite{Osip}. However, we will present a different proof which may be generalized to the vector-valued functions $f$. Since all functions we consider are defined on $[0,1]$, we will omit this fact in our notation.

We denote by $\eps_s$ a sequence of Rademacher functions (they should not be confused with $r_k$; this is another copy of a sequence of Rademacher functions and me may assume that they are defined on some other probability space $\Omega$). If $X$ is a Banach space and $f$ is an $X$-valued function, then the analogue of inequality \eqref{RdF_walsh} has the following form:
\begin{equation}
	\Big\| \sum_s \eps_s P_{I_s} f \Big\|_{L^p (\Rad X)}\lesssim \|f\|_{L^p(X)}.
	\label{RdF_walsh_Banach}
\end{equation}
We denote by $\Rad X$ the closure in $L^p(\Omega; X)$ of $X$-valued functions of the form 
$$
\sum_{j=1}^k  \eps_j(\omega) x_j, \quad x_j\in X.
$$
It is not difficult to see that Khintchine--Kahane inequality (see for instance \cite[p. 191]{HytBook}) implies that this definition does not depend on $p$ for $1\le p<\infty$. Also, for $X=\mathbb{R}$ Khintchine's inequality implies that the inequalities \eqref{RdF_walsh} and \eqref{RdF_walsh_Banach} are equivalent.

We introduce the following definition.

\begin{defin}
	We say that Banach space $X$ has the $\LPR_p^w$ property if the inequality \eqref{RdF_walsh_Banach} holds for any function $f\in L^p(X)$.
\end{defin}

This is an analogue of the $\LPR_p$ property, which was axiomatised in \cite{BGT},  for Walsh functions. The $\LPR_p$ property (which deals with the Fourier transform and the original inequality from the paper \cite{RdF} in the vector-valued setting) was studied in the papers \cite{HTY}, \cite{PSX} and others. In this paper we prove the results of the paper \cite{PSX} for our setting of Walsh functions. Now we pass to the exact statements of our main results.

\begin{theorem}
	If $X$ is a Banach lattice such that its 2-concavification $X_{(2)}$ is a UMD Banach lattice then $X$ is a space with $\LPR_p^w$ property for $2<p<\infty$.
	\label{thm1}
\end{theorem}

\begin{theorem}
	If $X$ is a Banach space with $\LPR_p^w$ property for some $p\ge 2$, then $X$ is also a Banach space with $\LPR_q^w$ property for any $q>p$.
\end{theorem}

These are the analogues of the results of the paper \cite{PSX} for $\LPR_p$ property. All necessary definitions regarding Banach lattices may be found in the book \cite{LinTzaf}. 




\section{$LPR_p^w$ property in Banach lattices}

In this section we are going to prove Theorem \ref{thm1}. We note that for UMD Banach lattices the inequality \eqref{RdF_walsh_Banach} can be rewritten in the same form as the inequality for scalar-valued functions:
$$
	\Big\|\Big(\sum_s |P_{I_s}f|^2\Big)^{1/2}\Big\|_{L^p(X)}\lesssim \|f\|_{L^p(X)}.
$$
This fact follows from the Khintchine--Maurey inequality (see for  instance \cite[p.86]{HytBook2}; note that UMD property implies finite cotype).

We start with the presentation of the proof of the desired inequality for scalar-valued functions which is different from the proof in the paper \cite{Osip} and then generalize it to the case of $X$-valued functions $f$.

\subsection{A combinatorial construction}

Suppose that $I_s=[a_s, b_s)$. We are going to use the decomposition of these intervals which was constructed in the paper \cite{Osip}:
\begin{equation}
I_s=\{a_s\}\cup \bigcup_{j\in\Theta_s}J_{js}\cup\bigcup_{i\in \widetilde{\Theta}_s} \wt{J}_{is}, \quad \text{where}\quad \Theta_s, \wt{\Theta}_s\subset\mathbb{N}.
\label{comb_const}
\end{equation}
The intervals $J_{js}$ and $\wt{J}_{is}$ are such that the following relations hold:
\begin{equation}
a_s\dotplus J_{js}=\delta_j \quad \text{and} \quad b_s\dotplus\wt{J}_{is}=\delta_i.
\label{decomp}
\end{equation}
We briefly describe the construction of such decomposition (and refer the reader to the paper \cite{Osip} for details).

Let us omit the index $s$ and construct the decomposition of an interval $I=[a,b)$. Consider the binary expansion of the numbers $a$ and $b$:
$$
a=\sum_{k=0}^N \alpha_{k}2^k, \quad b=\sum_{k=0}^N \beta_{k}2^k.
$$
First, we decompose the interval $[0, b)$. We number all digits in the binary expansion of $b$ which are equal to $1$:
$$
\beta_{k_1}=\beta_{k_2}=\ldots=\beta_{k_l}=1, \quad k_1>k_2>\ldots>k_l.
$$
Now we take the following intervals:
$$
\wt{J}_{k_i+1}=\Big[ \sum_{r=1}^{i-1}2^{k_r},\  \sum_{t=1}^{i} 2^{k_t} \Big), \quad i=1, 2, \ldots, l.
$$
It is not difficult to see that $b\dotplus \wt{J}_{k_i+1}=\delta_{k_i+1}$.

One of these segments contains the number $a$. Suppose that $a\in\wt{J}_{k_m+1}$. It means that $\alpha_j=\beta_j$ for $j>k_m$ and $\alpha_{k_m}=0$. Now we need to decompose the interval
$$
\Big(a,\ \sum_{t=1}^{m} 2^{k_t}\Big).
$$
It may be done in a similar way. We number the binary digits of $a$ which are equal to $0$:
$$
\alpha_{k_m}=\alpha_{\kap_h}=\ldots=\alpha_{\kap_1}=0, \quad k_m>\kap_h>\ldots>\kap_1.
$$
Now we take the following intervals:
$$
J_{\kap_i+1}=\Big[ a+\sum_{j=1}^{i-1} 2^{\kap_j} + 1, \  a+\sum_{g=0}^{i} 2^{\kap_g}  \Big], \quad i=1, 2,  \ldots, h.
$$
It is not difficult to see that $a\dotplus J_{\kap_i+1}=\delta_{\kap_i+1}$ and therefore we get the desired decomposition.

\subsection{The main argument for scalar-valued functions}

The inequality \eqref{RdF_walsh} will now follow from the following three inequalities:
\begin{align*}
\Big\|\Big(\sum_s |P_{\{a_s\}}f|^2\Big)^{1/2}\Big\|_{L^p}&\lesssim \|f\|_{L^p},\\
\Big\|\Big(\sum_s \Big|\sum_{j\in\Theta_s}P_{J_{js}}f\Big|^2\Big)^{1/2}\Big\|_{L^p}&\lesssim \|f\|_{L^p}, \\
 \Big\|\Big(\sum_s \Big|\sum_{i\in\wt{\Theta}_s}P_{\wt{J}_{is}}f\Big|^2\Big)^{1/2}\Big\|_{L^p}&\lesssim \|f\|_{L^p}.
\end{align*}
The first inequality follows easily from orthogonality arguments (because $\sum_s |P_{\{a_s\}} f|^2\le \|f\|_{L^2}^2$). Therefore, we should prove the second inequality (since the third one has the similar form and can be proved in the same way). We introduce the following notation: 
$$
J_s=\bigcup_{j\in\Theta_s} J_{js}.
$$
We have:
$$
P_{J_s}f=w_{a_s}\sum_{j\in\Theta_s} \Delta_j(w_{a_s}f).
$$
This identity is an easy consequence of the formulas \eqref{decomp}.
Therefore, we can rewrite our inequality in the following form:
$$
\Big\| \Big( \sum_s \Big| \sum_{j\in\Theta_s} \Delta_j (w_{a_s}f) \Big|^2 \Big)^{1/2} \Big\|_{L^p}\lesssim \|f\|_{L^p}.
$$

We denote by $G$ the following operator from $L^p$ to the space of $\ell^2$-valued functions:
\begin{equation}
f\mapsto \Big( \sum_{j\in\Theta_s} \Delta_j(w_{a_s}f) \Big)_s.
\label{Gf_def}
\end{equation}
We need to show that $G$ is a bounded operator from $L^p$ to $L^p(\ell^2)$. For $\ell^2$-valued function $g$ define the (dyadic) sharp maximal function in a usual way:
$$
g^\#(x)=\sup_{I\ni x} \Big( \frac{1}{|I|} \int_{I} \|g(t)-g_I\|_{\ell^2}^2\, dt \Big)^{1/2},\quad \text{where} \ g_{I}=\frac{1}{|I|}\int_{I}g(s)\, ds \in\ell^2.
$$
The supremum here is taken over all dyadic intervals which contain the point $x$.

We also need the following dyadic maximal function (at this time it is defined for scalar-valued functions $f$):
$$
M_2 f(x)=\sup_{I\ni x} \Big(\frac{1}{|I|}\int_I |f|^2\Big)^{1/2}.
$$
Here the supremum is also taken over all dyadic intervals which contain the point $x$.

Our goal is to prove the following pointwise estimate:
\begin{equation}
(Gf)^\# (x)\lesssim M_2f(x).
\label{key_est}
\end{equation}
This estimate will finish the proof of the theorem for scalar-valued functions $f$ because we may write:
$$
\|Gf\|_{L^p(\ell^2)}\lesssim \|(Gf)^\#\|_{L^p}\lesssim \|M_2 f\|_{L^p}\lesssim \|f\|_{L^p}.
$$
The first inequality here follows from the fact that the inequality
$$
\|g\|_{L^p(\ell^2)}\lesssim \|g^\#\|_{L^p}
$$
holds for any $\ell^2$-valued function $g$ such that $\int_0^1 g = 0.$ Note that $\int_0^1 Gf=0\in\ell^2$ because $\int_0^1 \Delta_j(h)=0$ for any function $h$ and any number $j>0$.

In order to prove the inequality \eqref{key_est}, it is enough to show that the following estimate holds: 
\begin{equation}
\Big( \frac{1}{|I|}\int_I \sum_s \Big| \sum_{j\in\Theta_s} \Delta_j(w_{a_s}f)-h_s \Big|^2 \Big)^{1/2}\lesssim \Big(\frac{1}{|I|}\int_I |f|^2\Big)^{1/2}.
\label{Gf-sharp}
\end{equation}
where
$$
h_s=\frac{1}{|I|}\int_I \sum_{j\in\Theta_s}\Delta_j(w_{a_s}f).
$$
Here $I$ is a dyadic interval.

We note that for any function $g$ the integral $\int_I \Delta_j(g)$ is equal to $0$ whenever $2^{j-1}\ge |I|^{-1}$ --- this is true because for any such interval $\int_I \E_j g=\int_I \E_{j-1} g=\int_I g$. Therefore, the expression for $h_s$ can be rewritten in the following way:
$$
h_s=\frac{1}{|I|}\int_I \sum_{\substack{j\in\Theta_s,\\ 2^j\le |I|^{-1}}}\Delta_j(w_{a_s}f).
$$
On the other hand, for any function $g$ the function $\Delta_j (g)$ is constant on the interval $I$ whenever $|I|\le 2^{-j}$ (because both functions $\E_j g$ and $\E_{j-1}g$ are constants on $I$). Hence, the left-hand side in the formula \eqref{Gf-sharp} equals to the following expression:
$$
\Big( \frac{1}{|I|}\int_I \sum_s \Big| \sum_{\substack{j\in\Theta_s, \\2^{j-1}\ge |I|^{-1}}} \Delta_j(w_{a_s}f)\Big|^2 \Big)^{1/2}. 
$$
We consider the restriction of the function $f$ to the interval $I$: $\wt{f}=f|_I$. Suppose that $|I|=2^{-m}$ and
$$
a_s=2^{k_1}+2^{k_2}+\ldots+2^{k_t}
$$
is the dyadic expansion of one of the numbers $a_s$, where
$$ 
 k_1<\ldots < k_t, \quad k_l\le m-1,\ k_{l+1}\ge m.
$$
We denote by $\wt{a}_s$ the number
$$
2^{k_{l+1}-m}+\ldots+2^{k_t-m}.
$$
Also, we denote by $\wt{w}_n$ the Walsh functions ``scaled" to the interval $I$ (i.e., they form the orthonormal basis in the space $L^2(I; |I|^{-1}dx)$). Furthermore, we use the notation $\wt{\Delta}_j$ for the dyadic martingale differences on $I$. Then the function $\Delta_j(w_{a_s} f)$ coincides on $I$ with the function $\pm\wt{\Delta}_{j-m}(\wt{w}_{\wt{a}_s}\wt{f})$ (note that the functions $r_{k_1+1}, \ldots, r_{k_l+1}$ are constants equal to $\pm 1$ on $I$; also, note that the operators $\Delta_j$ are ``local" in a sense that the value of a function $\Delta_j g$ on $I$ depends only on the restriction of $g$ to the interval $I$ if $2^{j-1}\ge |I|$). 

Now it is not difficult to see that the fact that 
$$
\{a_s\dotplus\delta_j\}_{j\in\Theta_s}
$$
are pairwise disjoint sets for different values of $s$ and $j$ (which is true because $a_s\dotplus\delta_j=J_{js}$; see the formula \eqref{decomp}) implies that the sets
$$
\{\wt{a}_s\dotplus \delta_{j-m}\}_{\substack{j\in\Theta_s,\\j\ge m+1}}
$$
are also pairwise disjoint.

Therefore, simply using Plancherel theorem, wee see that
$$
\Big( \frac{1}{|I|}\int_I \sum_s \Big| \sum_{\substack{j\in\Theta_s, \\2^{j-1}\ge |I|^{-1}}} \Delta_j(w_{a_s}f)\Big|^2 \Big)^{1/2}\leq \|\wt{f}\|_{L^2(I; |I^{-1}|dx)}=\Big(\frac{1}{|I|}\int_I |f|^2\Big)^{1/2},
$$
and the proof of the inequality \eqref{key_est} is finished.

\subsection{The inequality for vector-valued functions}

Once the estimate \eqref{key_est} is proved, we may complete the proof of Theorem \ref{thm1} in a similar way as it is done in the paper \cite{PSX} (for the usual $\LPR_p$ property).

Suppose that $X$ is a Banach lattice. We may assume that $X$ is a K\"othe function space defined on some probability space $(\Omega, \Sigma, \mu)$ (see \cite[p.25]{LinTzaf}). Indeed, the UMD property for the lattice $X_{(2)}$ implies the UMD property for the lattice $X$ (see \cite[Theorem 4]{RdF2}) and hence $X$ is reflexive (see \cite[page 183]{Pis}); besides that, we may assume that $X$ is separable and that it contains a weak unit.

Recall that $2$-concavification of $X$ is a Banach lattice $X_{(2)}$ with the norm 
$$
\|x\|_{X_{(2)}}=\||x|^{1/2}\|_X^2.
$$
The space $X_{(2)}$ is a Banach lattice if and only if $X$ is $2$-convex, that is, the following inequality holds for any $x_1,\ldots, x_n\in X$:
$$
\Big\| \Big( \sum_{j=1}^n |x_j|^2 \Big)^{1/2} \Big\|_X \le \Big( \sum_{j=1}^n \|x_j\|_X^2 \Big)^{1/2}.
$$

It is worth noting that the overview \cite{RdF2} is an excellent reference for many facts about Banach lattices with UMD property which we are going to use.

As we obtained before (see the beginning of Subsection 2.2), Theorem 1 would follow once we prove the following three inequalities (this time for $X$-valued functions $f$):
\begin{align}
	\Big\|\Big(\sum_s |P_{\{a_s\}}f|^2\Big)^{1/2}\Big\|_{L^p(X)}&\lesssim \|f\|_{L^p(X)},\label{easy_part}\\
	\Big\|\Big(\sum_s \Big|\sum_{j\in\Theta_s}P_{J_{js}}f\Big|^2\Big)^{1/2}\Big\|_{L^p(X)}&\lesssim \|f\|_{L^p(X)}, \label{main_part}\\
	\Big\|\Big(\sum_s \Big|\sum_{i\in\wt{\Theta}_s}P_{\wt{J}_{is}}f\Big|^2\Big)^{1/2}\Big\|_{L^p(X)}&\lesssim \|f\|_{L^p(X)}.\label{analog_part}
\end{align}

The estimate \eqref{easy_part} is again easy. Note that for any fixed $\omega\in\Omega$ Plancherel theorem implies the inequality
$$
\Big( \sum_s |(f,w_{a_s})|^2 \Big)^{1/2}(\omega)\lesssim \Big( \int_0^1 |f|^2 \Big)^{1/2}(\omega).
$$
Therefore, we may write:
\begin{multline*}
\Big\|	\Big( \sum_s |(f,w_{a_s})|^2 \Big)^{1/2}\Big\|_X\lesssim \Big\| \Big( \int_0^1 |f|^2 \Big)^{1/2} \Big\|_X\leq \Big( \int_0^1 \|f\|_X^2 \Big)^{1/2}\\ \leq \Big( \int_0^1 \|f\|_X^p \Big)^{1/p}=\|f\|_{L^p(X)}.
\end{multline*}
We used $2$-convexity of $X$ in the second inequality (note that we may assume that $f$ is a finite linear combination of Walsh functions and therefore it is constant on the dyadic intervals of length $2^{-n}$ for sufficiently large value of $n$; hence, we can write the sum instead of the integral and use $2$-convexity).

Now we pass to the proof of the estimates \eqref{main_part} and \eqref{analog_part}. We only show how to prove the inequality \eqref{main_part} since the proof of \eqref{analog_part} is similar.

Recall that the space $X(\ell^2)$ is defined as the space of sequences $x=(x_j)\subset X$ such that
$$
\|x\|_{X(\ell^2)}=\sup_{n\in\mathbb{N}}\Big\|\Big(\sum_{j=1}^n |x_j|^2\Big)^{1/2}\Big\|_X<\infty.
$$
We define the operator $G$ in the same way as we did for scalar-valued functions (see the formula \eqref{Gf_def}). This time, our goal is to prove that $G$ is a bounded operator from $L^p(X)$ to $L^p(X(\ell^2))$.

For an $X$-valued function $f$ we define the maximal function $M_2f$ in the following way:
$$
M_2 f(x)=\sup_{I\ni x} \Big(\frac{1}{|I|}\int_I |f|^2\Big)^{1/2}.
$$
Note that this is a function of two variables: for each fixed $\omega\in\Omega$ this is a maximal function of $f(\cdot, \omega)$.

For an $X(\ell^2)$-valued function $g$ we define the dyadic sharp maximal function $g^\#$ in a similar way: for each fixed $\omega\in\Omega$ we apply the usual sharp maximal function (which we already defined in the previous Subsection) to the $\ell^2$-valued function $g(\cdot, \omega)$. The formula which defines $g^\#$ for $X(\ell^2)$-valued function $g=(g_1, g_2, \ldots)$ has the following form:
$$
g^\#(x)=\sup_{I\ni x} \Big( \frac{1}{|I|} \int_{I} \sum_{j=1}^\infty |g_j(t)-(g_I)_j|^2\, dt \Big)^{1/2},\quad \text{where} \ g_{I}=\frac{1}{|I|}\int_{I}g(s)\, ds \in X(\ell^2).
$$
The inequality \eqref{key_est} for scalar valued functions implies the following estimate in our context:
$$
G(f(\cdot, \omega))^\#\lesssim M_2(f(\cdot, \omega)), \quad \text{a.e.}\ \omega\in\Omega.
$$
Therefore, we have the following estimate:
$$
\|(Gf)^\#\|_{L^p(X)}\lesssim \|M_2 f\|_{L^p(X)}.
$$
Hence, the proof of Theorem 1 is finished once we prove the following two estimates:
\begin{equation}
\|g\|_{L^p(X(\ell^2))}\lesssim \|g^\#\|_{L^p(X)}\quad \text{and} \quad \|M_2 f\|_{L^p(X)}\lesssim \|f\|_{L^p(X)}
\label{stand_est}
\end{equation}
for an arbitrary $X(\ell^2)$-valued function $g$ such that $\int g=0$ (recall that in the previous subsection we proved that $\int Gf=0$) and $X$-valued function $f$. The second inequality follows immediately from the boundedness of Hardy--Littlewood maximal operator in $L^p(X)$ for UMD Banach lattices $X$ (which is proved in \cite{BourMaks}. See also \cite[Theorem 3]{RdF2}), we simply need to apply it to $X_{(2)}$ here:
$$
\|M_2f\|_{L^p(X)}^2=\|M[|f|^2]\|_{L^{p/2}(X_{(2)})}\lesssim \||f|^2\|_{L^{p/2}(X_{2})}=\|f\|_{L^p(X)}^2.
$$
Here $Mf$ is the standard (dyadic) Hardy--Littlewood maximal function: 
$$
Mf(x)=\sup_{I\ni x} \frac{1}{|I|}\int_I |f(s)|\, ds.
$$
Note that $M_2f(x)=M[|f^2|]^{1/2}$.

Now we explain how to obtain the first of the inequalities \eqref{stand_est} for an $X(\ell^2)$-valued function $g=(g_1, g_2, \ldots) $. We use the following standard inequality for $\ell^2$-valued functions $u$ and $v$ such that $\int_0^1 v = 0$ (see for instance \cite[p. 223]{MusSchlag} for scalar-valued functions; as usual, there is no difference with $\ell^2$-valued case): 
$$
\int_0^1 |\langle u(x), v(x) \rangle|\, dx \lesssim \int_0^1 (\Sq u)(x)v^\#(x)\, dx.
$$
Here $\Sq$ denotes the martingale square function operator with respect to the Haar filtration.

We take any function $h\in L^{p'}(X^*(\ell^2))$ and write:
\begin{equation}
	\label{final}
\Big|\int_{[0,1]\times\Omega} \langle g, h\rangle\Big|\lesssim \int_{[0,1]\times\Omega} g^\# \Sq h\lesssim \|g^\#\|_{L^p(X)}\|\Sq h\|_{L^{p'}(X^*)}\lesssim \|g^\#\|_{L^p(X)}\|h\|_{L^{p'}(X^*(\ell^2))}.
\end{equation}
In the last inequality we used the following simple statement (we note that if $X$ is a UMD Banach lattice then $X^*$ is also a Banach lattice with UMD property).

\begin{lemma*}
	If $Y$ is a UMD Banach lattice and $h=(h_1, h_2, \ldots)$ is a $Y(\ell^2)$-valued function then for any $q$, $1<q<\infty$, we have the following inequality: $\|\Sq h\|_{L^q(Y)}\lesssim \|h\|_{L^q(Y(\ell^2))}$.
\end{lemma*}

Since it turned out to be difficult to find an exact reference for this statement, we present its proof here.

\begin{proof}
	We again suppose that $Y$ is a Banach function space defined on some measure space $(\Omega', \Sigma', \mu')$. Denote by $d_j$ the operators of martingale difference with respect to Haar filtration. Using Khintchine's inequality and then Minkowski's integral inequality we can write the following estimate:
	\begin{multline*}
		\Sq h(x, \omega')=\Big( \sum_n\sum_j |d_j h_n(x, \omega')|^2 \Big)^{1/2}\lesssim \Big( \sum_n \Big[ \E \Big| \sum_j \eps_j d_j h_n(x, \omega') \Big| \Big]^2 \Big)^{1/2}\\
		\leq \E \Big( \Big[ \sum_n \Big| \sum_j \eps_j d_j h_n(x,\omega') \Big|^2 \Big]^{1/2} \Big).
	\end{multline*}

Using this inequality, we estimate the norm of the square function in the following way:
\begin{multline*}
	\|\Sq h\|_{L^q(Y)}^q \lesssim \int \Big\|\E \Big( \Big[ \sum_n \Big| \sum_j \eps_j d_j h_n(x) \Big|^2 \Big]^{1/2} \Big)\Big\|_Y^q dx\\
	\leq \int\E \Big[ \Big\| \Big( \sum_n \Big| \sum_j \eps_j d_j h_n(x) \Big|^2 \Big)^{1/2} \Big\|_Y^q \Big]\, dx=\E \Big( \Big\| \sum_j \eps_j d_j h \Big\|_{L^q(Y(\ell^2))}^q \Big)\lesssim \|h\|_{L^q(Y(\ell^2))}^q.
\end{multline*}
Here the last inequality follows from the definition of UMD property for the space $Y(\ell^2)$ (which follows from the UMD property of $Y$, see the overview \cite{RdF2}).

\end{proof}

Now we take the supremum in the estimate \eqref{final} over all  $h\in L^{p'}(X^*(\ell^2))$ with the unit norm and get the desired bound.

\section{$LPR_p^w$ implies $LPR_q^w$ for $q>p$}

Now we pass to the proof of Theorem 2. Suppose that $X$ is a Banach space with $\LPR_p^w$ property (recall that it means that the inequality \eqref{RdF_walsh_Banach} holds for any $f\in L^p(X)$). We note that it implies that $X$ is a UMD space (since we may take the ``family" $\{I_s\}$ consisting of one interval $[0,N]$ and therefore get the uniform boundedness of projections $P_{[0,N]}$ which already implies UMD property; see for instance \cite[p.254]{Pis}).

Our goal is to prove that for any family of disjoint intervals $\{I_s\}$ the following inequality holds for $f\in L^q(X)$:
$$
\Big\| \sum_s \eps_s P_{I_s} f \Big\|_{L^q (\Rad X)}\lesssim \|f\|_{L^q(X)}, \ q>p.
$$
We again use the combinatorial construction from the paper \cite{Osip} (see \eqref{comb_const} and \eqref{decomp}). We introduce the notation $J_{0s}=a_s$. It is enough to prove the following two inequalities:
\begin{align*}
	\Big\| \sum_s \eps_s \sum_{j\in\Theta_s\cup\{0\}} P_{J_{js}}f \Big\|_{L^q(\Rad X)}\lesssim \|f\|_{L^q(X)},\\
	\Big\| \sum_s \eps_s \sum_{i\in\wt{\Theta}_s} P_{\wt{J}_{is}} f \Big\|_{L^q(\Rad X)}\lesssim \|f\|_{L^q(X)}.
\end{align*}

Since the proofs of these inequalities are similar, we only show how to prove the first of them. Recall that $P_{J_{js}}f=w_{a_s}\Delta_j(w_{a_s}f)$ and therefore  our inequality may be rewritten in the following form:
$$
	\Big\| \sum_s \eps_s w_{a_s} \sum_{j\in\Theta_s\cup\{0\}} \Delta_j(w_{a_s}f) \Big\|_{L^q(\Rad X)}\lesssim \|f\|_{L^q(X)}.
$$
Now we use the standard Kahane's contraction principle (see \cite[p.181]{HytBook}) and conclude that our inequality is equivalent to the following:
$$
	\Big\| \sum_s \eps_s  \sum_{j\in\Theta_s\cup\{0\}} \Delta_j(w_{a_s}f) \Big\|_{L^q(\Rad X)}\lesssim \|f\|_{L^q(X)}.
$$
Consider the operator $T$ which maps the function $f\in L^p(X)$ to
$$
\sum_s \eps_s  \sum_{j\in\Theta_s\cup\{0\}} \Delta_j(w_{a_s}f).
$$
We see that $\LPR_p^w$ property implies that $T$ is a bounded linear operator from $L^p(X)$ to $L^p(\Rad X)$ (here it is important to note that the set $\cup_{j\in{\Theta_s\cup\{0\}}} J_{js}$ is a segment in $\Z_+$) and we need to prove that $T$ is bounded from $L^q(X)$ to $L^q(\Rad X)$. 

We consider the operator $T^*: L^{p'}(\Rad X^*)\to L^{p'}(X^*)$ (note that $(\Rad X)^* \sim \Rad X^*$, see Section 3 in \cite{HytWeis}). The direct computation shows that $T^*$ has the following form:
\begin{equation}
T^*\Big( \sum_s \eps_s g_s \Big)=\sum_{s}w_{a_s} \sum_{j\in\Theta_s\cup\{0\}}\Delta_j g_s.
\label{form_of_T^*}
\end{equation}
Here $g=\sum_s\eps_s g_s$ is a $\Rad X^*$-valued function. In order to prove that $T^*$ is a bounded operator from $L^{q'}(\Rad X^*)$ to $L^{q'}(X^*)$ we use the following vector-valued version of Calder\'{o}n--Zygmund decomposition from the paper \cite{Kis}.

\begin{fact}
	Suppose that $E$ is a Banach space, $g$ is a simple $E$-valued function (which means that $g=\E_k g$ for sufficiently large values of $k$) and $\lambda$ is an arbitrary positive number. Then there exist simple functions $b$ and $h$ such that $g=b+h$ and the following conditions hold.
	\begin{enumerate}
		\item $\|h\|_{L^\infty(Y)}\lesssim \lambda$ and $\|h\|_{L^1(Y)}\lesssim \|g\|_{L^1(Y)}$.
		
		\item $\int_0^1 b = 0$ and for every $n\ge 1$ we have $\Delta_n b = \mathbbm{1}_{e_n}\Delta_n b$, where $e_n\in\F_{n-1}$ and $|\bigcup_{n\ge 1} e_n|\lesssim \lambda^{-1}\|g\|_{L^1(Y)}$.
	\end{enumerate}
\end{fact}

Our goal is to prove that the operator $T^*$ is of weak type $(1,1)$ because in this case the boundedness from $L^{q'}(\Rad X^*)$ to $L^{q'}(X^*)$ would follow by Marcinkiewicz interpolation theorem. We fix a number $\lambda>0$ and apply Fact 1 to $Y=\Rad X^*$ and our function $g$. Note that we may consider only simple functions $g$ since by density we may assume that the collection $\{I_s\}$ is finite and each function $g_s$ is a finite linear combination of Walsh functions. We have:
\begin{multline}
|\{x: \|(T^*g)(x)\|_{X^*} > \lambda\}| \\ \leq |\{x: \|(T^*b)(x)\|_{X^*} > \lambda/2\}| + |\{x: \|(T^*h)(x)\|_{X^*} > \lambda/2\}|.
\label{ff}
\end{multline}
In order to estimate the second summand in this formula, we use the boundedness of $T^*$ from $L^{p'}(\Rad X^*)$ to $L^{p'}(X^*)$:
\begin{multline*}
|\{x: \|(T^*h)(x)\|_{X^*} > \lambda/2\}|\leq \Big(\frac{\lambda}{2}\Big)^{-p'}\|T^* h\|_{L^{p'}(X^*)}^{p'} \lesssim \lambda^{-p'} \|h\|_{L^{p'}(\Rad X^*)}^{p'}\\ \leq \lambda^{-p'}\|h\|_{L^\infty(\Rad X^*)}^{p'-1}\|h\|_{L^1(\Rad X^*)}\lesssim \lambda^{-1}\|g\|_{L^1(\Rad X^*)}. 
\end{multline*}

Now we estimate the first summand from the formula \eqref{ff}. It is easy to see that the formula \eqref{form_of_T^*} and the property of $b$ from Fact 1 imply that the supports of the function $T^*b$ is contained in the set $\bigcup_{n\ge 1} e_n$. Therefore, we have:
$$
|\{x: \|(T^*b)(x)\|_{X^*} > \lambda/2\}|\leq |\{x: (T^*b)(x)\neq 0\}|\leq \Big| \bigcup_{n\ge 1} e_n \Big|\lesssim \lambda^{-1} \|g\|_{L^1(\Rad X^*)},
$$
and the proof of Theorem 2 is finished.

\bigskip

St. Petersburg Department, Steklov Math. Institute, Fontanka 27, St. Petersburg 191023 Russia

\medskip

celis-anton@yandex.ru

\end{document}